\documentclass[12pt]{article}%
\usepackage{amsfonts}
\usepackage{sw20bams}
\usepackage{amsmath}
\usepackage{amssymb}
\usepackage{graphicx}%
\providecommand{\U}[1]{\protect\rule{.1in}{.1in}}
\begin{document}

\title{Fractional Iterates and Oscillatory Convergence}
\author{Steven Finch}
\date{October 29, 2026}
\maketitle

\begin{abstract}
The simple continued fractions for the Golden \&\ Silver means are well-known.
\ It is astonishing that, as far as we know, no one has published
half-iterates (let alone quarter-iterates) for the corresponding algorithms.
We also examine the cosine and logistic maps (with parameter $2<\lambda<3$).\ 

\end{abstract}

\footnotetext{Copyright \copyright \ 2025-2026 by Steven R. Finch. All rights
reserved.}

For each of the iterative examples we have examined in the past
\cite{F1-oscill, F2-oscill, F3-oscill, F4-oscill, F5-oscill}, the convergence
to a fixed point has been monotonic. \ A question -- what happens in the event
that convergence is oscillatory? -- remains unanswered. \ Our purpose here is
to extend Koenig's method \cite{Ax-oscill} so that real (not
complex)\ fractional iterates can be constructed. \ As far as is known, this
approach is new. [Corrigendum at end.]

\section{$1+1/x$}

\ \ \ \ We begin with standard techniques yielding complex results, before
transitioning to a nonstandard procedure yielding real results.

\subsection{Schr\"{o}der}

Consider the recurrence
\[%
\begin{array}
[c]{ccccc}%
x_{k}=f(x_{k-1})=1+\dfrac{1}{x_{k-1}} &  & \text{for }k\geq1\text{;} &  &
x_{0}=1.
\end{array}
\]
The function $f(x)$ has an attracting fixed point at $x=\varphi=\left(
1+\sqrt{5}\right)  /2$, the Golden mean. \ Set $y=x-\varphi$ and note
\[%
\begin{array}
[c]{ccccc}%
y_{k}+\varphi=1+\dfrac{1}{y_{k-1}+\varphi} &  & \text{hence define} &  &
g(y)=1-\varphi+\dfrac{1}{y+\varphi}.
\end{array}
\]
We have $-1<g^{\prime}(0)=-1/(1+\varphi)<0$. The solution of Schr\"{o}der's
equation \cite{Ax-oscill}
\[
G\left(  g(y)\right)  =g^{\prime}(0)G\left(  y\right)
\]
is obtained via Koenig's method either iteratively:%
\[
G\left(  y\right)  =\lim\limits_{k\rightarrow\infty}\dfrac{y_{k}}{g^{\prime
}(0)^{k}}=\lim\limits_{k\rightarrow\infty}\,(-1)^{k}(1+\varphi)^{k}y_{k}%
\]
or via power series (matching coefficients and summing a geometric series):%
\[
G(y)=\frac{\sqrt{5}\,y}{\sqrt{5}+y}%
\]
where $y=y_{0}$. \ Backtracking:%
\[%
\begin{array}
[c]{ccc}%
K(x)=G(x-\varphi)=\dfrac{\sqrt{5}\,\left(  x-\varphi\right)  }{\sqrt
{5}+\left(  x-\varphi\right)  }, &  & K^{[-1]}(x)=\dfrac{-\varphi\,x+\sqrt
{5}\,\left(  x+\varphi\right)  }{\sqrt{5}-x}%
\end{array}
\]
it follows that%
\[
x_{1/2}=f^{[1/2]}(1)=K^{[-1]}\left(  \left(  \dfrac{-1}{1+\varphi}\right)
^{1/2}K(1)\right)  =\frac{3}{2}-\frac{1}{2}i,
\]%
\[
x_{3/2}=f^{[3/2]}(1)=K^{[-1]}\left(  \left(  \dfrac{-1}{1+\varphi}\right)
^{3/2}K(1)\right)  =\frac{8}{5}+\frac{1}{5}i
\]
where $i$ is the imaginary unit. \ Such formulas were anticipated by Helms
\cite{Hm-oscill}, appearing instead as%
\[%
\begin{array}
[c]{ccc}%
\dfrac{1}{x_{1/2}}=\dfrac{3}{5}+\dfrac{1}{5}i, &  & \dfrac{1}{x_{3/2}}%
=\dfrac{8}{13}-\dfrac{1}{13}i.
\end{array}
\]
More generally, starting from arbitrary $x_{0}=x$,
\[
x_{1/2}=f^{[1/2]}(x)=K^{[-1]}\left(  \left(  \dfrac{-1}{1+\varphi}\right)
^{1/2}K(x)\right)  =\frac{1+(1+i)x}{i+x},
\]%
\[
x_{3/2}=f^{[3/2]}(x)=K^{[-1]}\left(  \left(  \dfrac{-1}{1+\varphi}\right)
^{3/2}K(x)\right)  =\frac{(1+i)+(2+i)x}{1+(1+i)x}%
\]
which enable quick tests that $f^{[1/2]}(f^{[1/2]}(x))=f(x)$ and
$f^{[1/2]}(f^{[3/2]}(x))=f^{[2]}(x)$.

\subsection{Abel}

Let $D$ denote the disk in the complex plane centered at point $1/2$ and of
radius $\varphi-1/2$. \ The boundary $\partial D$ of $D$ is the circle passing
through four points $1-\varphi$, $1/2\pm(\varphi-1/2)i$ and $\varphi$. \ Let%
\[
\log_{\theta}(z)=\ln\left\vert z\right\vert +i\arg_{\theta}(z)
\]
be the complex logarithm with branch cut at angle $\theta$. \ Define%
\[
F_{\theta}(z)=\frac{\log_{\theta}\left(  \dfrac{\sqrt{5}\,\left(
z-\varphi\right)  }{\sqrt{5}+\left(  z-\varphi\right)  }\right)  }%
{\log_{\theta}\left(  \dfrac{-1}{1+\varphi}\right)  }.
\]
One solution of Abel's equation \cite{Ax-oscill}%
\[
F\left(  1+\dfrac{1}{z}\right)  =F(z)+1
\]
is $F_{+\pi/2}(z)$, valid for $\left\vert z-1/2\right\vert <\varphi-1/2$;
\ another is $F_{-\pi/2}(z)$, valid for $\left\vert z-1/2\right\vert
>\varphi-1/2$. \ These facts can be experimentally verified via computer
algebra. \ To clarify: $\theta$ is the lower bound of the polar-angle interval
of length $2\pi$, e.g.,

\begin{itemize}
\item $f(1/2)=3$ and $\arg_{+\pi/2}(1/2)=\arg_{+\pi/2}(3)=2\pi,$

\item $f(-1/2)=-1$ and $\arg_{+\pi/2}(-1/2)=\arg_{+\pi/2}(-1)=\pi,$

\item $f(2)=3/2$ and $\arg_{-\pi/2}(2)=\arg_{-\pi/2}(3/2)=0,$

\item $f(-2)=1/2$ and $\arg_{-\pi/2}(-2)=\pi>\arg_{-\pi/2}(1/2)=0.$
\end{itemize}

\noindent Obviously the verification fails for $z=\varphi$ \& $z=1-\varphi$
and perhaps too for certain other $z\in\partial D$. \ Let $z_{0}=1$
\&\ $z_{1}=2$. \ We have%
\[%
\begin{array}
[c]{ccc}%
F\left(  \dfrac{3}{2}-\dfrac{1}{2}i,\dfrac{\pi}{2}\right)  =F\left(
1,\dfrac{\pi}{2}\right)  +\dfrac{1}{2}, &  & F\left(  \dfrac{8}{5}+\dfrac
{1}{5}i,\dfrac{\pi}{2}\right)  =F\left(  1,\dfrac{\pi}{2}\right)  +\dfrac
{3}{2}%
\end{array}
\]
implying that Schr\"{o}der-based values $x_{1/2}$ \& $x_{3/2}$ coincide with
Abel-based values $z_{1/2}$ \&\ $z_{3/2}$ when $\theta=\pi/2$. \ Further,
$z_{1/2},z_{3/2}\in\partial D$ and $F(z_{1/2},\pi/2)=F(z_{1/2},-\pi/2)$, but%
\[
F\left(  z_{3/2},\dfrac{\pi}{2}\right)  -F\left(  z_{3/2},-\dfrac{\pi}%
{2}\right)  =\frac{2\pi}{\pi+i\ln(1+\varphi)}\neq0.
\]
Details of complex iterates as such are left for readers. \ We move on to real iterates.

\subsection{Modified Abel}

This is the genesis of our excursion into unfamiliar territory. \ From
Fibonacci ratios:%
\[
x_{0}=1<\dfrac{3}{2}=x_{2}<\ldots<\varphi<\ldots<x_{3}=\dfrac{5}{3}<2=x_{1}%
\]
it is natural to define%
\[%
\begin{array}
[c]{ccc}%
g(x)=(\varphi-1)-\dfrac{\varphi-x}{\varphi+1-x}, &  & h(x)=-(\varphi
-1)+\dfrac{\varphi+x}{\varphi+1+x}%
\end{array}
\]
which capture the separate behaviors of $x_{k}$ for even subscripts and odd
subscripts, respectively \cite{F6-oscill}. \ We have $0<g^{\prime
}(0)=h^{\prime}(0)=1/(1+\varphi)^{2}<1$. \ Koenig's method yields%
\[%
\begin{array}
[c]{ccc}%
G(x)=\dfrac{\sqrt{5}\,x}{\sqrt{5}-x}, &  & H(x)=\dfrac{\sqrt{5}\,x}{\sqrt
{5}+x}.
\end{array}
\]
Let%
\[
F(x)=\left\{
\begin{array}
[c]{lll}%
\dfrac{\kappa+\ln\left(  \dfrac{\sqrt{5}\,\left(  \varphi-x\right)  }{\sqrt
{5}-\left(  \varphi-x\right)  }\right)  -\ln\left(  \dfrac{\sqrt{5}}{\sqrt
{5}-1}\right)  }{\ln\left(  \dfrac{1}{(1+\varphi)^{2}}\right)  } &  & \text{if
}x<\varphi,\smallskip\\
\dfrac{\ln\left(  \dfrac{\sqrt{5}\,\left(  x-\varphi\right)  }{\sqrt
{5}+\left(  x-\varphi\right)  }\right)  -\ln\left(  \dfrac{\sqrt{5}}{\sqrt
{5}+1}\right)  }{\ln\left(  \dfrac{1}{(1+\varphi)^{2}}\right)  } &  & \text{if
}x>\varphi
\end{array}
\right.
\]
where $\kappa$ is a parameter. \ For $f(x)=1+1/x$, the equation%
\[
F(f(x))=F(x)+\Delta
\]
possesses no solution unless $\Delta=1/2$ and $\kappa=\ln(1+\varphi)$. \ To
determine half-iterates for oscillatory convergence, the usual additive term
$1/2$ is therefore replaced by $1/4$. \ It follows that%
\[
x_{0}=1<x_{1/2}=f^{[1/2]}(1)=F^{[-1]}\left(  F(1)+\frac{1}{4}\right)
=\frac{4}{2}-\dfrac{1}{2}\varphi=1.19098...,
\]%
\[
x_{1}=2>x_{3/2}=f^{[3/2]}(1)=F^{[-1]}\left(  F(2)+\frac{1}{4}\right)
=\dfrac{17}{11}+\dfrac{2}{11}\varphi=1.83964...,
\]%
\[
x_{2}=\frac{3}{2}<x_{5/2}=f^{[5/2]}(1)=F^{[-1]}\left(  F\left(  \frac{3}%
{2}\right)  +\frac{1}{4}\right)  =\dfrac{48}{29}-\dfrac{2}{29}\varphi
=1.54358...,
\]%
\[
x_{3}=\frac{5}{3}>x_{7/2}=f^{[7/2]}(1)=F^{[-1]}\left(  F\left(  \frac{5}%
{3}\right)  +\frac{1}{4}\right)  =\dfrac{61}{38}+\dfrac{1}{38}\varphi
=1.64784...,
\]%
\[
x_{4}=\frac{8}{5}<x_{9/2}=f^{[9/2]}(1)=F^{[-1]}\left(  F\left(  \frac{8}%
{5}\right)  +\frac{1}{4}\right)  =\dfrac{323}{199}-\dfrac{2}{199}%
\varphi=1.60685...
\]
and $x_{5}=13/8$. \ Likewise,%
\begin{align*}
x_{0}  &  =1<x_{1/3}=f^{[1/3]}(1)=F^{[-1]}\left(  F(1)+\frac{1}{6}\right) \\
&  =\frac{1}{7}\left[  11-(20+35\varphi)^{1/3}+(-55+35\varphi)^{1/3}\right]
=1.13280...,
\end{align*}%
\[
x_{2/3}=F^{[-1]}\left(  F(x_{1/3})+\frac{1}{6}\right)  =\frac{1}{5}\left[
8-(3+5\varphi)^{1/3}+(-8+5\varphi)^{1/3}\right]  =1.24367...
\]
and%
\begin{align*}
x_{0}  &  =1<x_{1/4}=f^{[1/4]}(1)=F^{[-1]}\left(  F(1)+\frac{1}{8}\right) \\
&  =\frac{1}{11}\left[  17+2\varphi-\sqrt{-15+50\varphi}\right]  =1.10164...,
\end{align*}%
\[
x_{3/4}=F^{[-1]}\left(  F(x_{1/2})+\frac{1}{8}\right)  =\frac{1}{29}\left[
46+2\varphi-\sqrt{-55+130\varphi}\right]  =1.26801....
\]
The function $F^{[-1]}$ actually has two values ($F(1)=F(3)$, $F(2)=F(4/3)$,
etc.), one given above and the other given below. \ A parallel set of results
exist for Lucas ratios:%
\[
\tilde{x}_{0}=3>\tilde{x}_{1/2}=f^{[1/2]}(3)=F^{[-1]}\left(  F(3)+\frac{1}%
{4}\right)  =\dfrac{3}{2}+\dfrac{1}{2}\varphi=2.30901...,
\]%
\[
\tilde{x}_{1}=\frac{4}{3}<\tilde{x}_{3/2}=f^{[3/2]}(3)=F^{[-1]}\left(
F\left(  \frac{4}{3}\right)  +\frac{1}{4}\right)  =\dfrac{19}{11}-\dfrac
{2}{11}\varphi=1.43308...,
\]%
\[
\tilde{x}_{2}=\frac{7}{4}>\tilde{x}_{5/2}=f^{[5/2]}(3)=F^{[-1]}\left(
F\left(  \frac{7}{4}\right)  +\frac{1}{4}\right)  =\dfrac{46}{29}+\dfrac
{2}{29}\varphi=1.69779...,
\]%
\[
\tilde{x}_{3}=\frac{11}{7}<\tilde{x}_{7/2}=f^{[7/2]}(3)=F^{[-1]}\left(
F\left(  \frac{11}{7}\right)  +\frac{1}{4}\right)  =\dfrac{62}{38}-\dfrac
{1}{38}\varphi=1.58899...,
\]%
\[
\tilde{x}_{4}=\frac{18}{11}>\tilde{x}_{9/2}=f^{[9/2]}(3)=F^{[-1]}\left(
F\left(  \frac{18}{11}\right)  +\frac{1}{4}\right)  =\dfrac{321}{199}%
+\dfrac{2}{199}\varphi=1.62932...
\]
and $\tilde{x}_{5}=29/18$. \ Likewise,%
\begin{align*}
\tilde{x}_{0}  &  =3>\tilde{x}_{1/3}=f^{[1/3]}(3)=F^{[-1]}\left(
F(3)+\frac{1}{6}\right) \\
&  =\frac{1}{3}\left[  5+(2+3\varphi)^{1/3}+(5-3\varphi)^{1/3}\right]
=2.47532...,
\end{align*}%
\[
\tilde{x}_{2/3}=F^{[-1]}\left(  F(\tilde{x}_{1/3})+\frac{1}{6}\right)
=\frac{1}{11}\left[  18+(35+55\varphi)^{1/3}+(90-55\varphi)^{1/3}\right]
=2.18083...
\]
and%
\begin{align*}
\tilde{x}_{0}  &  =3>\tilde{x}_{1/4}=f^{[1/4]}(3)=F^{[-1]}\left(
F(3)+\frac{1}{8}\right) \\
&  =\frac{1}{11}\left[  17+2\varphi+\sqrt{-15+50\varphi}\right]  =2.57764...,
\end{align*}%
\[
\tilde{x}_{3/4}=F^{[-1]}\left(  F(\tilde{x}_{1/2})+\frac{1}{8}\right)
=\frac{1}{29}\left[  46+2\varphi+\sqrt{-55+130\varphi}\right]  =2.12757....
\]
A plot of $x$ \&\ $\tilde{x}$ points suggests that fractional interpolation
here consists of disjoint curved segments from $[k,k+1)$ for each integer
$k\geq0$. \ The large discontinuities at $k=1$ \&\ $k=2$ are especially unsatisfying.

\section{$2+1/x$}

Our focus henceforth will be on the modified Abel approach. \ The function
$f(x)=2+1/x$ has an attracting fixed point at $x=\psi=1+\sqrt{2}$, the Silver
mean. \ From Pell ratios:%
\[
x_{0}=2<\dfrac{12}{5}=x_{2}<\ldots<\psi<\ldots<x_{3}=\dfrac{29}{12}<\dfrac
{5}{2}=x_{1}%
\]
it is natural to define%
\[%
\begin{array}
[c]{ccc}%
g(x)=(\psi-2)-\dfrac{\psi-x}{2\psi+1-2x}, &  & h(x)=\dfrac{\psi+x}{2\psi
+1+2x}-(\psi-2)
\end{array}
\]
which capture the separate behaviors of $x_{k}$ for even subscripts and odd
subscripts, respectively \cite{F6-oscill}. \ We have $0<g^{\prime
}(0)=h^{\prime}(0)=1/(1+2\psi)^{2}<1$. \ Koenig's method yields
\[%
\begin{array}
[c]{ccc}%
G(x)=\dfrac{2\sqrt{2}\,x}{2\sqrt{2}-x}, &  & H(x)=\dfrac{2\sqrt{2}\,x}%
{2\sqrt{2}+x}.
\end{array}
\]
Let%
\[
F(x)=\left\{
\begin{array}
[c]{lll}%
\dfrac{\kappa+\ln\left(  \dfrac{2\sqrt{2}\,\left(  \psi-x\right)  }{2\sqrt
{2}-\left(  \psi-x\right)  }\right)  -\ln\left(  \dfrac{2\sqrt{2}}{2\sqrt
{2}-1}\right)  }{\ln\left(  \dfrac{1}{(1+2\psi)^{2}}\right)  } &  & \text{if
}x<\psi,\smallskip\\
\dfrac{\ln\left(  \dfrac{2\sqrt{2}\,\left(  x-\psi\right)  }{2\sqrt{2}+\left(
x-\psi\right)  }\right)  -\ln\left(  \dfrac{2\sqrt{2}}{2\sqrt{2}+1}\right)
}{\ln\left(  \dfrac{1}{(1+2\psi)^{2}}\right)  } &  & \text{if }x>\psi
\end{array}
\right.
\]
where $\kappa$ is a parameter. \ The equation
\[
F(f(x))=F(x)+\Delta
\]
possesses no solution unless $\Delta=1/2$ and $\kappa=\ln(5/7+(4/7)\psi)$.
\ To determine half-iterates for oscillatory convergence, the usual additive
term $1/2$ is therefore replaced by $1/4$. \ It follows that%
\[
x_{0}=2<x_{1/2}=f^{[1/2]}(2)=F^{[-1]}\left(  F(2)+\frac{1}{4}\right)
=\dfrac{18}{7}-\dfrac{1}{7}\psi=2.22654...,
\]%
\[
x_{1}=\frac{5}{2}>x_{3/2}=f^{[3/2]}(2)=F^{[-1]}\left(  F\left(  \frac{5}%
{2}\right)  +\frac{1}{4}\right)  =\dfrac{98}{41}+\dfrac{1}{41}\psi
=2.44912...,
\]%
\[
x_{2}=\frac{12}{5}<x_{5/2}=f^{[5/2]}(2)=F^{[-1]}\left(  F\left(  \frac{12}%
{5}\right)  +\frac{1}{4}\right)  =\dfrac{578}{239}-\dfrac{1}{239}%
\psi=2.40830...,
\]%
\[
x_{3}=\frac{29}{12}>x_{7/2}=f^{[7/2]}(2)=F^{[-1]}\left(  F\left(  \frac
{29}{12}\right)  +\frac{1}{4}\right)  =\dfrac{3362}{1393}+\dfrac{1}{1393}%
\psi=2.41522...,
\]%
\[
x_{4}=\frac{70}{29}<x_{9/2}=f^{[9/2]}(2)=F^{[-1]}\left(  F\left(  \frac
{70}{29}\right)  +\frac{1}{4}\right)  =\dfrac{19602}{8119}-\dfrac{1}{8119}%
\psi=2.41403...
\]
and $x_{5}=169/70$. \ Likewise,%
\begin{align*}
x_{0}  &  =2<x_{1/3}=f^{[1/3]}(2)=F^{[-1]}\left(  F(2)+\frac{1}{6}\right) \\
&  =\frac{1}{17}\left[  41-(14+34\psi)^{1/3}+(-82+34\psi)^{1/3}\right]
=2.16802...,
\end{align*}%
\[
x_{2/3}=F^{[-1]}\left(  F(x_{1/3})+\frac{1}{6}\right)  =\frac{1}{29}\left[
70-(12+29\psi)^{1/3}+(-70+29\psi)^{1/3}\right]  =2.27191...
\]
and%
\begin{align*}
x_{0}  &  =2<x_{1/4}=f^{[1/4]}(2)=F^{[-1]}\left(  F(2)+\frac{1}{8}\right) \\
&  =\frac{1}{41}\left[  98+\psi-2\sqrt{-28+29\psi}\right]  =2.13294...,
\end{align*}%
\[
x_{3/4}=F^{[-1]}\left(  F(x_{1/2})+\frac{1}{8}\right)  =\frac{1}{239}\left[
576+\psi-2\sqrt{-168+169\psi}\right]  =2.29050....
\]
The function $F^{[-1]}$ actually has two values ($F(2)=F(3)$, $F(5/2)=F(7/3)$,
etc.), one given above and the other given below. \ A parallel set of results
exist for Pell-Lucas ratios:%
\[
\tilde{x}_{0}=3>\tilde{x}_{1/2}=f^{[1/2]}(3)=F^{[-1]}\left(  F(3)+\frac{1}%
{4}\right)  =\dfrac{16}{7}+\dfrac{1}{7}\psi=2.63060...,
\]%
\[
\tilde{x}_{1}=\frac{7}{3}<\tilde{x}_{3/2}=f^{[3/2]}(3)=F^{[-1]}\left(
F\left(  \frac{7}{3}\right)  +\frac{1}{4}\right)  =\dfrac{100}{41}-\dfrac
{1}{41}\psi=2.38014..,
\]%
\[
\tilde{x}_{2}=\frac{17}{7}>\tilde{x}_{5/2}=f^{[5/2]}(3)=F^{[-1]}\left(
F\left(  \frac{17}{7}\right)  +\frac{1}{4}\right)  =\dfrac{576}{239}+\dfrac
{1}{239}\psi=2.42014...,
\]%
\[
\tilde{x}_{3}=\frac{41}{17}<\tilde{x}_{7/2}=f^{[7/2]}(3)=F^{[-1]}\left(
F\left(  \frac{41}{17}\right)  +\frac{1}{4}\right)  =\dfrac{3364}{1393}%
-\dfrac{1}{1393}\psi=2.41319...,
\]%
\[
\tilde{x}_{4}=\frac{99}{41}>\tilde{x}_{9/2}=f^{[9/2]}(3)=F^{[-1]}\left(
F\left(  \frac{99}{41}\right)  +\frac{1}{4}\right)  =\dfrac{19600}%
{8119}+\dfrac{1}{8119}\psi=2.41438...
\]
and $\tilde{x}_{5}=239/99$. \ Likewise,%
\begin{align*}
\tilde{x}_{0}  &  =3>\tilde{x}_{1/3}=f^{[1/3]}(3)=F^{[-1]}\left(
F(3)+\frac{1}{6}\right) \\
&  =\frac{1}{36}\left[  87+(783-324\psi)^{1/3}+(135+324\psi)^{1/3}\right]
=2.71228...,
\end{align*}%
\[
\tilde{x}_{2/3}=F^{[-1]}\left(  F(\tilde{x}_{1/3})+\frac{1}{6}\right)
=\frac{1}{41}\left[  99+(198-82\psi)^{1/3}+(34+82\psi)^{1/3}\right]
=2.57243...
\]
and%
\begin{align*}
\tilde{x}_{0}  &  =3>\tilde{x}_{1/4}=f^{[1/4]}(3)=F^{[-1]}\left(
F(3)+\frac{1}{8}\right) \\
&  =\frac{1}{41}\left[  98+\psi+2\sqrt{-28+29\psi}\right]  =2.76530...,
\end{align*}%
\[
\tilde{x}_{3/4}=F^{[-1]}\left(  F(\tilde{x}_{1/2})+\frac{1}{8}\right)
=\frac{1}{239}\left[  576+\psi+2\sqrt{-168+169\psi}\right]  =2.54978....
\]
Upon graphing, the sizeable breaks at integers $k$ between curved segments are
again disappointing. \ This may be an artifact of our fractional model $F$,
equipped with just two parameters $\{\kappa,\Delta\}$. Or it might be true
that a real continuous fractional interpolation is plainly impossible in this
case. \ 

\section{$\cos(x)$}

Some background is provided in \cite{MO-oscill}. \ From%
\[
x_{0}=0<0.54\approx\cos(1)=x_{2}<\ldots<\theta<\ldots<x_{3}=\cos
(\cos(1))\approx0.85<1=x_{1}%
\]
it is natural to define
\[%
\begin{array}
[c]{ccc}%
g(x)=\theta-\cos(\cos(\theta-x)), &  & h(x)=\cos(\cos(\theta+x))-\theta
\end{array}
\]
which capture the separate behaviors of $x_{k}$ for even subscripts and odd
subscripts, respectively \cite{F7-oscill}. \ The limiting value%
\[
\theta=0.7390851332151606416553120...
\]
is Dottie's number \cite{Kapl-oscill, Pain-oscill}. \ We have $0<g^{\prime
}(0)=h^{\prime}(0)=1-\theta^{2}<1$. \ Koenig's method yields%
\[%
\begin{array}
[c]{ccc}%
G(x)=x+%
{\displaystyle\sum\limits_{j=2}^{\infty}}
\gamma_{j}x^{j}, &  & H(x)=-G(-x)=x+%
{\displaystyle\sum\limits_{j=2}^{\infty}}
(-1)^{j-1}\gamma_{j}x^{j}%
\end{array}
\]
where%
\[
\gamma_{2}=\dfrac{1-\theta^{2}-\sqrt{1-\theta^{2}}}{2\theta\left(
1-\theta^{2}\right)  }=-0.3277931305953677271045803...,
\]%
\[
\gamma_{3}=\dfrac{2-2\theta^{2}-3\sqrt{1-\theta^{2}}}{6\theta^{2}\left(
1-\theta^{2}\right)  }=-0.7486243776658610242164076...,
\]%
\begin{align*}
\gamma_{4}  &  =\dfrac{22-33\theta^{2}+20\theta^{4}-6\theta^{6}-\left(
22-31\theta^{2}+11\theta^{4}\right)  \sqrt{1-\theta^{2}}}{24\theta^{3}\left(
1-\theta^{2}\right)  \left(  3-3\theta^{2}+\theta^{4}\right)  }\\
&  =0.4577246514588478179217120...,
\end{align*}%
\[%
\begin{array}
[c]{ccc}%
\gamma_{5}=0.3010324155830439550743236..., &  & \gamma_{6}%
=-0.4470215225416579495572650...,
\end{array}
\]%
\[%
\begin{array}
[c]{ccc}%
\gamma_{7}=-0.1066634338781197816943286..., &  & \gamma_{8}%
=0.3716684648125405104090988....
\end{array}
\]
The power series for $G(x)$ is slowly convergent:\ using 16 terms, only three
digits of \cite{F7-oscill}%
\[
G(\theta-0)=0.3983002403035094139563243...
\]
are correctly predicted. \ Doing the same for $H(x)$, nine digits of%
\[
H(1-\theta)=0.2682998330950090571338993...
\]
are correctly predicted (due to smallness of $1-\theta$ relative to $\theta$).
We must therefore abandon the matching-coefficients method and adopt
iteration, so as to implement the modified Abel approach with required accuracy.

\section{$\lambda\,x\,(1-x)$, $2<\lambda<3$}

For simplicity, let $\lambda=5/2$. \ From%
\[
x_{0}=\frac{1}{2}<0.59\approx\frac{75}{128}=x_{2}<\ldots<\mu<\ldots
<x_{3}=\frac{19875}{32768}\approx0.61<\frac{5}{8}=x_{1}%
\]
it is natural to define%
\[
g(x)=(\lambda-2)^{2}x-(\lambda-3)(\lambda-2)\lambda\,x^{2}-2(\lambda
-2)\lambda^{2}x^{3}-\lambda^{3}x^{4},
\]%
\[
h(x)=(\lambda-2)^{2}x+(\lambda-3)(\lambda-2)\lambda\,x^{2}-2(\lambda
-2)\lambda^{2}x^{3}+\lambda^{3}x^{4}%
\]
which capture the separate behaviors of $x_{k}$ for even subscripts and odd
subscripts, respectively \cite{F7-oscill}. \ The limiting value $\mu
=(\lambda-1)/\lambda$ is $3/5$ when $\lambda=5/2$. \ We have $0<g^{\prime
}(0)=h^{\prime}(0)=(\lambda-2)^{2}<1$. \ Koenig's method yields%
\[%
\begin{array}
[c]{ccc}%
G(x)=x+%
{\displaystyle\sum\limits_{j=2}^{\infty}}
\gamma_{j}x^{j}, &  & H(x)=-G(-x)=x+%
{\displaystyle\sum\limits_{j=2}^{\infty}}
(-1)^{j-1}\gamma_{j}x^{j}%
\end{array}
\]
where%
\[%
\begin{array}
[c]{ccccccc}%
\gamma_{2}=-\dfrac{10}{3}, &  & \gamma_{3}=-\dfrac{200}{9}, &  & \gamma
_{4}=\dfrac{1000}{9} &  & \gamma_{5}=-\dfrac{4000}{27}%
\end{array}
\]%
\[%
\begin{array}
[c]{ccccc}%
\gamma_{6}=-\dfrac{1000000}{891}, &  & \gamma_{7}=\dfrac{68000000}{18711}, &
& \gamma_{8}=\dfrac{4810000000}{344817}.
\end{array}
\]
The power series for $G(x)$ is slowly convergent, although faster than that
examined in Section 3:\ using 16 terms, six digits of \cite{F7-oscill}%
\[
G(\mu-x_{0})=G(1/10)=0.0533831106341909825926069...
\]
are correctly predicted. \ Doing the same for $H(x)$, seventeen digits of%
\[
H(x_{1}-\mu)=H(1/40)=0.0266915553170954912963034...
\]
are correctly predicted (due to smallness of $x_{1}-\mu$ relative to
$\mu-x_{0}$). \ A symbolic expression for the modified Abel solution $F(x)$,
as in Sections 1.3 \&\ 2, is an unrealistic hope. \ Coupled numerically with
either matching-coefficients or iterative methods, however, the approach would
seem to be feasible.

Apart from some preliminary analysis in \cite{F7-oscill}, the case $\lambda=3$
is wide open. \ Its limiting value is $2/3$ and associated odd/even
recurrences are%
\[%
\begin{array}
[c]{ccc}%
u_{k}=u_{k-1}-18u_{k-1}^{3}-27u_{k-1}^{4}, &  & u_{0}=1/12;
\end{array}
\]%
\[%
\begin{array}
[c]{ccc}%
v_{k}=v_{k-1}-18v_{k-1}^{3}+27v_{k-1}^{4}, &  & v_{0}=1/6.
\end{array}
\]
Such iterations are more challenging than others in the present paper (due to
the unit coefficients for $u_{k-1}$, $v_{k-1}$ and missing $u_{k-1}^{2}$,
$v_{k-1}^{2}$ terms). \ Relevant constants are%
\[%
\begin{array}
[c]{ccc}%
C_{u}=-0.1805303007686495535981970..., &  & C_{v}%
=-0.1388636341019828869315303...
\end{array}
\]
and we wonder about the algebraic independence of these.

\section{Corrigendum}

The following arose from simplifying the expression $F^{[-1]}(F(x)+1/4)$ in
Section 1.3:%
\[%
\begin{array}
[c]{ccc}%
p(x)=\dfrac{(\varphi-1)+2\,\varphi\,x}{(\varphi+1)+(\varphi-1)x}, &  &
q(x)=\dfrac{\varphi-2(\varphi-1)x}{(\varphi-2)+\varphi\,x}.
\end{array}
\]
These satisfy, for instance,%
\[%
\begin{array}
[c]{ccc}%
p\left(  x_{0}\right)  =p(1)=\dfrac{4}{2}-\dfrac{1}{2}\varphi=x_{1/2}, &  &
q\left(  \tilde{x}_{1/2}\right)  =q\left(  \dfrac{3}{2}+\dfrac{1}{2}%
\varphi\right)  =\dfrac{4}{3}=\tilde{x}_{1}%
\end{array}
\]
but%
\[%
\begin{array}
[c]{ccc}%
p\left(  x_{1/2}\right)  =\dfrac{4}{3}\neq2=x_{1}, &  & q\left(  \tilde{x}%
_{0}\right)  =q(3)=\dfrac{4}{2}-\dfrac{1}{2}\varphi\neq\dfrac{3}{2}+\dfrac
{1}{2}\varphi=\tilde{x}_{1/2}%
\end{array}
\]
are contrary to (flawed) expectation. \ More generally,%
\[
q(p(x))=1+\frac{1}{x}=p(q(x))
\]
but $p(p(x))$ \&\ $q(q(x))$ do not \ resemble $f(x)$ at all. \ We have
therefore \textit{not} found a real half-iterate for the Golden-mean
continued-fraction algorithm. \ The same regrettable conclusion applies to
Section 2:%
\[%
\begin{array}
[c]{ccc}%
r(x)=\dfrac{(-\psi+5)+(\psi+9)x}{(3\psi-1)+(-\psi+5)x}, &  & s(x)=\dfrac
{(2\psi-1)+7\,\psi\,x}{(3\psi+2)+(2\psi-1)x}%
\end{array}
\]
however for a slightly different reason. \ These satisfy, for instance,%
\[%
\begin{array}
[c]{ccc}%
r\left(  x_{0}\right)  =r(2)=\dfrac{18}{7}-\dfrac{1}{7}\psi=x_{1/2}, &  &
s\left(  \tilde{x}_{0}\right)  =s(3)=\dfrac{16}{7}+\dfrac{1}{7}\psi=\tilde
{x}_{1/2}%
\end{array}
\]
but%
\[%
\begin{array}
[c]{ccc}%
r\left(  x_{1/2}\right)  =\dfrac{7}{3}\neq\dfrac{5}{2}=x_{1}, &  & s\left(
\tilde{x}_{1/2}\right)  =\dfrac{5}{2}\neq\dfrac{7}{3}=\tilde{x}_{1}%
\end{array}
\]
are again contrary to earlier misconception. \ Here none of $s(r(x))$,
$r(s(x))$, $r(r(x))$,\ $s(s(x))$ resemble $f(x)$ at all. \ We have therefore
\textit{not} found a real half-iterate for the Silver-mean continued-fraction
algorithm. \ 

Certain formulas in Section 1.1 appeared, in fact, at least as far back as
1879. \ Given $f(x)=(a\,x+b)/(c\,x+d)$, Johnson \cite{John-oscill} obtained
two half-iterates:%
\[%
\begin{array}
[c]{ccc}%
f_{+}^{[1/2]}(x)=\dfrac{\left(  a+\sqrt{a\,d-b\,c}\right)  x+b}{c\,x+\left(
d+\sqrt{a\,d-b\,c}\right)  }, &  & f_{-}^{[1/2]}(x)=\dfrac{\left(
a-\sqrt{a\,d-b\,c}\right)  x+b}{c\,x+\left(  d-\sqrt{a\,d-b\,c}\right)  }%
\end{array}
\]
which yield%
\[%
\begin{array}
[c]{ccc}%
f_{\pm}^{[1/2]}(1)=\dfrac{3}{2}\mp\dfrac{1}{2}i, &  & f_{\pm}^{[3/2]}%
(1)=\dfrac{8}{5}\pm\dfrac{1}{5}i
\end{array}
\]
if $a=b=c=1$ \&\ $d=0$. \ Changing $a$ to $2$ yields $f_{\pm}^{[1/2]}(1)=2\mp
i$ instead. \ We nearly included $2-i$ in Section 2 but chose not to (in
fruitless pursuit of real iterates). \ \clearpage

\section{Addendum}

We prove that the matching-coefficient technique indeed gives a geometric
series. \ Define $\varphi=\left(  1+\sqrt{5}\right)  /2$ and $\xi=-1/\sqrt{5}%
$. \ \ From%
\[
\varphi+\xi=\frac{\sqrt{5}+5-2}{2\sqrt{5}}=\frac{3+\sqrt{5}}{2\sqrt{5}}%
=\frac{1+2\sqrt{5}+5}{4\sqrt{5}}=-\varphi^{2}\xi
\]
we obtain $\varphi^{-2}(\varphi+\xi)=-\xi$, hence $\varphi^{-2m}(\varphi
+\xi)^{m}=(-\xi)^{m}$ for all integers $m$. \ This lemma will be helpful
later. \ Let%
\[
g(x)=1-\varphi+\dfrac{1}{x+\varphi}%
\]
which has an an attracting fixed point at $x=0$. \ The solution of
Schr\"{o}der's functional equation%
\[
G\left(  g(x)\right)  =-\frac{1}{\varphi^{2}}G\left(  x\right)
\]
via power series is our focus. \ Let%
\[
g(x)=%
{\displaystyle\sum\limits_{n=1}^{\infty}}
a_{n}x^{n}=%
{\displaystyle\sum\limits_{n=1}^{\infty}}
\frac{(-1)^{n}}{\varphi^{n+1}}x^{n}%
\]
and, for arbitrary $c_{1}$,%
\[
G(y)=%
{\displaystyle\sum\limits_{k=1}^{\infty}}
c_{k}y^{k}.
\]
The $n^{\text{th}}$ coefficient of $g(x)^{k}$, where $n\geq k$, is%
\[
b_{n,k}=%
{\displaystyle\sum\limits_{\substack{m_{1}+m_{2}+\cdots+m_{k}=n,\\\text{each
}m_{j}\geq1}}}
a_{m_{1}}a_{m_{2}}\cdots a_{m_{k}}=\frac{(-1)^{n}}{\varphi^{n+k}}\dbinom
{n-1}{k-1}%
\]
because clearly
\[
a_{m_{1}}a_{m_{2}}\cdots a_{m_{k}}=\frac{(-1)^{m_{1}+m_{2}+\cdots+m_{k}}%
}{\varphi^{m_{1}+m_{2}+\cdots+m_{k}+k}}=\frac{(-1)^{n}}{\varphi^{n+k}}%
\]
and \cite{Aw-oscill} the number of integer compositions of $n$ into $k$
positive parts is $\tbinom{n-1}{k-1}$. \ The matching-coefficient condition
for each $n\geq1$ is thus%
\[%
{\displaystyle\sum\limits_{k=1}^{n}}
\frac{(-1)^{n}}{\varphi^{n+k}}\dbinom{n-1}{k-1}c_{k}=-\frac{1}{\varphi^{2}%
}c_{n},
\]
i.e.,%
\[
(-1)^{n}%
{\displaystyle\sum\limits_{k=1}^{n}}
\dbinom{n-1}{k-1}\varphi^{2-n-k}c_{k}=-c_{n}.
\]
For example, if $n=1$, then%
\[
(-1)^{1}\varphi^{2-1-1}c_{1}=-c_{1}%
\]
(which conveys no new information)\ and, if $n=2$,%
\[
\varphi^{2-2-1}c_{1}+\varphi^{2-2-2}c_{2}=-c_{2}%
\]
which by the lemma implies%
\[
-\xi\,c_{1}-\varphi^{-2}\xi\,c_{1}+\varphi^{-2}c_{2}=\varphi^{-2}\left[
(\varphi+\xi)-\xi\right]  c_{1}+\varphi^{-2}c_{2}=-c_{2}%
\]
i.e., $c_{2}=\xi\,c_{1}$. \ As another example, if $n=3$, then%
\[
-\varphi^{2-3-1}c_{1}-2\varphi^{2-3-2}c_{2}-\varphi^{2-3-3}c_{3}=-c_{3}%
\]
which by the lemma implies%
\[
\xi^{2}c_{1}-\varphi^{-4}\xi^{2}c_{1}+\varphi^{-4}c_{3}=\varphi^{-4}\left[
(\varphi+\xi)^{2}-\xi^{2}\right]  c_{1}+\varphi^{-4}c_{3}=c_{3}%
\]
i.e., $c_{3}=\xi^{2}c_{1}$. \ As yet another example, if $n=4$, then%
\[
\varphi^{2-4-1}c_{1}+3\varphi^{2-4-2}c_{2}+3\varphi^{2-4-3}c_{3}%
+\varphi^{2-4-4}c_{4}=-c_{4}%
\]
which by the lemma implies
\[
-\xi^{3}c_{1}-\varphi^{-6}\xi^{3}c_{1}+\varphi^{-6}c_{4}=\varphi^{-6}\left[
(\varphi+\xi)^{3}-\xi^{3}\right]  c_{1}+\varphi^{-6}c_{4}=-c_{4}%
\]
i.e., $c_{4}=\xi^{3}c_{1}$. \ 

Assume $c_{k}=\xi^{k-1}c_{1}$ for all $1\leq k\leq n-1$; we shall prove that
$c_{n}=\xi^{n-1}c_{1}$. \ Observe%
\begin{align*}%
{\displaystyle\sum\limits_{k=1}^{n-1}}
\dbinom{n-1}{k-1}\varphi^{2-n-k}\xi^{k-1}  &  =%
{\displaystyle\sum\limits_{j=0}^{n-2}}
\dbinom{n-1}{j}\varphi^{1-n-j}\xi^{j}\\
&  =\varphi^{2-2n}%
{\displaystyle\sum\limits_{j=0}^{n-2}}
\dbinom{n-1}{j}\varphi^{n-1-j}\xi^{j}\\
&  =\varphi^{2-2n}\left[
{\displaystyle\sum\limits_{j=0}^{n-1}}
\dbinom{n-1}{j}\varphi^{n-1-j}\xi^{j}-\xi^{n-1}\right] \\
&  =\varphi^{-2(n-1)}\left[  (\varphi+\xi)^{n-1}-\xi^{n-1}\right] \\
&  =(-\xi)^{n-1}-\varphi^{2-2n}\xi^{n-1}%
\end{align*}
by the binomial theorem and lemma; this yields%
\begin{align*}
c_{n}  &  =(-1)^{n-1}%
{\displaystyle\sum\limits_{k=1}^{n}}
\dbinom{n-1}{k-1}\varphi^{2-n-k}c_{k}\\
&  =(-1)^{n-1}\left[  (-\xi)^{n-1}c_{1}-\varphi^{2-2n}\xi^{n-1}c_{1}%
+\varphi^{2-2n}c_{n}\right] \\
&  =\xi^{n-1}c_{1}+(-1)^{n}\left[  \varphi^{2-2n}\xi^{n-1}c_{1}-\varphi
^{2-2n}c_{n}\right]
\end{align*}
hence%
\[
\left[  1+(-1)^{n}\varphi^{2-2n}\right]  c_{n}=\xi^{n-1}\left[  1+(-1)^{n}%
\varphi^{2-2n}\right]  c_{1}%
\]
hence $c_{n}=\xi^{n-1}c_{1}$. \ The formula follows for all $n$ by induction
and therefore%
\[
G(x)=%
{\displaystyle\sum\limits_{n=1}^{\infty}}
\xi^{n-1}c_{1}x^{n}=\frac{c_{1}\sqrt{5}\,x}{\sqrt{5}+x}.
\]

\section{Acknowledgements}

I\ am grateful to Witold Jarczyk \& Karol Baron \cite{BJ-oscill} for helpful
discussions. \ The creators of Mathematica earn my gratitude every day:\ this
paper could not have otherwise been written. \ It constitutes almost surely my
last words on the subject.


\begin{thebibliography}{99}                                                                                               %


\bibitem {F1-oscill}S. R. Finch, Compositional square roots of $\exp(x)$ and
$1+x^{2}$, arXiv:2504.19999.

\bibitem {F2-oscill}S. R. Finch, Half-iterates of $x(1+x)$, $\sin(x)$ and
$\exp(x/e)$, arXiv:2506.07625v1.

\bibitem {F3-oscill}S. R. Finch, Half-iterates of $x\exp(x)$, $x+1/x$ and
$\operatorname{arcsinh}(x)$, \newline arXiv:2506.07625v2.

\bibitem {F4-oscill}S. R. Finch, Half-iterates and delta conjectures, arXiv:2506.07625.

\bibitem {F5-oscill}S. R. Finch, Exercises in iterational asymptotics IV, arXiv:2509.24918.

\bibitem {Ax-oscill}D. S. Alexander, \textit{A History of Complex Dynamics.
From Schr\"{o}der to Fatou and Julia}, Friedr. Vieweg \& Sohn, 1994, pp.
46--49; MR1260930.

\bibitem {Hm-oscill}G. Helms, Fractional iteration of the function
$f(x)=1/(1+x)$, https://go.helms-net.de/math/tetdocs/FracIterAltGeom.htm.

\bibitem {F6-oscill}S. R. Finch, Exercises in iterational asymptotics III, arXiv:2503.13378.

\bibitem {MO-oscill}G. Edgar, S. Ivanov, J. D. Hamkins and other contributors,
How to solve $f(f(x))=\cos(x)$? https://mathoverflow.net/questions/17605/how-to-solve-ffx-cosx.

\bibitem {F7-oscill}S. R. Finch, Exercises in iterational asymptotics II, arXiv:2501.06065.

\bibitem {Kapl-oscill}S. R. Kaplan, The Dottie number, \textit{Math. Mag.} 80
(2007) 73--74.

\bibitem {Pain-oscill}J.-C. Pain, An exact series expansion for the Dottie
number, arXiv:2303.17962.

\bibitem {John-oscill}W. W. Johnson, Symbolic powers and roots of functions in
the form $(a\,x+b)/(c\,x+d)$, \textit{Messenger of Math.} 9 (1879) 99--103.

\bibitem {Aw-oscill}G. E. Andrews, \textit{The Theory of Partitions},
Addison-Wesley, 1976, pp. 54--55; MR0557013.

\bibitem {BJ-oscill}K. Baron and W. Jarczyk, Recent results on functional
equations in a single variable, perspectives and open problems,
\textit{Aequationes Math.} 61 (2001) 1--48; MR1820808.%

\begin{tabular}
[c]{lll}
& Steven Finch & \\
& MIT Sloan School of Management & \\
& Cambridge, MA, USA & \\
& \textit{steven\_finch\_math@outlook.com} &
\end{tabular}

\end{thebibliography}
\end{document}